\documentclass[12pt,twoside]{amsart}
\usepackage{amsmath, amsthm, amscd, amsfonts, amssymb, graphicx}
\usepackage[bookmarksnumbered, plainpages]{hyperref}

\newtheorem{thm}{Theorem}[section]

\newtheorem{defn}[thm]{Definition}

\numberwithin{equation}{section}

\begin{document}

\title{\bf Minimal translation surfaces with respect to semi-symmetric connections in $\mathbb{R}^3$ and $\mathbb{R}^3_1$}
\author{Yong Wang}

\thanks{{\scriptsize
\hskip -0.4 true cm \textit{2010 Mathematics Subject Classification:}
53C40; 53C42.
}}

\maketitle

\begin{abstract}
 In this paper, we define and classify minimal translation surfaces with respect to a kind of semi-symmetric metric connections and a kind of semi-symmetric non-metric connections in $\mathbb{R}^3$ and $\mathbb{R}^3_1$.
\end{abstract}

\vskip 0.2 true cm


\pagestyle{myheadings}
\markboth{\rightline {\scriptsize Wang}}
         {\leftline{\scriptsize Minimal translation surfaces with respect to semi-symmetric connections }}

\bigskip
\bigskip


\section{ Introduction}
\indent Minimal surfaces are among the most natural objects in differential geometry, and have been studied during the last two and half centuries since J. L.
Lagrange. In particular, minimal surfaces have encountered striking applications in other fields, like mathematical physics, conformal geometry, computer aided design, among others. In order to search for more minimal surfaces, some natural geometric assumptions arise.
Translation surfaces were studied in the Euclidean $3$-dimensional space and they are represented as graphs $z=f(x)+g(y)$, where
$f$ and $g$ are smooth functions. Scherk proved in 1835 that, besides the planes, the only minimal translation surfaces are the surfaces given by
$$z=\frac{1}{a}{\rm log}\mid\frac{{\rm cos}(ax)}{{\rm cos}(ay)}\mid,$$
where $a$ is a non-zero constant. Since then, minimal translation surfaces were generalized in several directions. For example, the Euclidean space
$\mathbb{R}^3$ was replaced with other spaces of dimension $3$- usually being $3$-dimensional Lie groups and the notion of translation was often replaced by
using the group operation. See for example \cite{ILM},\cite{Lop},\cite{LM},\cite{Seo},\cite{Y13} or \cite{YLK}. Another generalizations of Scherk surfaces are: affine translation surfaces in Euclidean $3$-space \cite{LY}, affine translation surfaces in affine $3$-dimensional space \cite{YF} and translation surfaces in Galilean $3$-space \cite{Y17}. On the other hand, Scherk surfaces were generalized to minimal translation surfaces in Euclidean spaces of arbitrary dimensions. See for example \cite{DVZ},\cite{MM}.\\
\indent H. A. Hayden introduced the notion of a semi-symmetric metric connection on a
Riemannian manifold \cite{HA}. K. Yano studied a Riemannian manifold endowed with
a semi-symmetric metric connection \cite{Ya}. Some properties of a Riemannian manifold
and a hypersurface of a Riemannian manifold with a semi-symmetric metric
connection were studied by T. Imai \cite{I1,I2}. Z. Nakao \cite{NA} studied submanifolds of
a Riemannian manifold with semi-symmetric metric connections. In \cite{GE},  Gozutok and Esin studied the tangent bundle of a hypersurface with semi-symmetric
metric connections. In \cite{De}, Demirbag investigated the properties of a weakly Ricci symmetric manifold admitting a semi-symmetric metric connection.
N. S. Agashe and
M. R. Chafle introduced the notion of a semisymmetric non-metric connection and
studied some of its properties and submanifolds of a Riemannian manifold with a
semi-symmetric non-metric connection \cite{AC1,AC2}. \\
 \indent In this paper, we define and classify minimal translation surfaces with respect to a kind of semi-symmetric metric connections and a kind of semi-symmetric non-metric connections in $\mathbb{R}^3$ and $\mathbb{R}^3_1$.
In Section 2,  we define and classify minimal translation surfaces with respect to a kind of semi-symmetric metric connections and a kind of semi-symmetric non-metric connections in $\mathbb{R}^3$.
In Section 3, we define and classify minimal translation surfaces with respect
to a kind of semi-symmetric metric connections and a kind of semi-symmetric non-metric connections in $\mathbb{R}^3_1$.

\vskip 1 true cm

\section{Minimal translation surfaces with respect to semi-symmetric connections in $\mathbb{R}^3$ }
\indent Let $\mathbb{R}^3$ be the $3$-dimensional Euclidean space with the canonical Euclidean metric $\widetilde{g}$. Let $X_1=\frac{\partial}{\partial x}$, $X_2=\frac{\partial}{\partial y}$, $X_3=\frac{\partial}{\partial z}$. Let $\nabla^L_{X_i}X_j=0$ for $1\leq i,j\leq 3$ be the Levi-Civita connection on
$\mathbb{R}^3$. We define a special semi-symmetric metric connection by
\begin{equation}
\nabla_{X}Y=\nabla^L_{X}Y+\widetilde{g}(Y,X_3)X-\widetilde{g}(X,Y)X_3.
\end{equation}
The connection $\nabla$ of $\mathbb{R}^3$ is given by
\begin{align}
&\nabla_{X_1}X_1=-X_3,~~\nabla_{X_1}X_2=0,~~\nabla_{X_1}X_3=X_1,\\\notag
&\nabla_{X_2}X_1=0,~~\nabla_{X_2}X_2=-X_3,~~\nabla_{X_2}X_3=X_2,\\\notag
&\nabla_{X_3}X_1=0,~~\nabla_{X_3}X_2=0,~~\nabla_{X_3}X_3=0.
\notag
\end{align}
\begin{defn}
A surface $M$ in $\mathbb{R}^3$ is a translation surface if it is given by an isometric immersion $\mathcal{F}: U\subset \mathbb{R}^2\rightarrow \mathbb{R}^3$ of the form
\begin{equation}
\mathcal{F}(u,v)=(u,v,f(u)+g(v)),~~~~~~~~({\rm Type }~I)
\end{equation}
or
\begin{equation}
\mathcal{F}(u,v)=(u,f(u)+g(v),v),~~~~~~~~({\rm Type }~II)
\end{equation}
or
\begin{equation}
\mathcal{F}(u,v)=(f(u)+g(v),u,v),~~~~~~~~({\rm Type }~III)
\end{equation}
where $f$ and $g$ are smooth functions on open sets of $\mathbb{R}$.
\end{defn}
Let $E_1=\mathcal{F}_u$, $E_2=\mathcal{F}_v$ and $\{E_1,E_2\}$ be the basis of $TM$ and $N$ be the unit normal vector field of $TM$ in $\mathbb{R}^3$. Let $p:T\mathbb{R}^3\mid_M\rightarrow TM$ be the projection. The we have the Gauss formula with respect to $\nabla$
\begin{equation}
\nabla_XY=p\nabla_XY+\sigma(X,Y)N,
\end{equation}
where $X,Y\in TM$ and $\sigma(X,Y)$ is the second fundamental form with respect to $\nabla$. In general, $\sigma(X,Y)\neq \sigma(Y,X)$.
Let $e_1$, $e_2$ be orthonormal basis of $TM$. We define the mean curvature of $M$ with respect to $\nabla$:
$H=\frac{1}{2}[\sigma(e_1,e_1)+\sigma(e_2,e_2)].$ Let
\begin{equation}
E=\widetilde{g}(\mathcal{F}_u,\mathcal{F}_u),~~F=\widetilde{g}(\mathcal{F}_u,\mathcal{F}_v),~~G=\widetilde{g}(\mathcal{F}_v,\mathcal{F}_v).
\end{equation}
Then
\begin{equation}
H=\frac{G\widetilde{g}(\nabla_{E_1}E_1,N)-F\widetilde{g}(\nabla_{E_1}E_2,N)-F\widetilde{g}(\nabla_{E_2}E_1,N)+E\widetilde{g}(\nabla_{E_2}E_2,N)}
{2(EG-F^2)}.
\end{equation}
We called that $M$ is minimal with respect to $\nabla$ if $H=0$. So by (2.8), $M$ is minimal with respect to $\nabla$ if and only if
\begin{equation}
{G\widetilde{g}(\nabla_{E_1}E_1,N)-F\widetilde{g}(\nabla_{E_1}E_2,N)-F\widetilde{g}(\nabla_{E_2}E_1,N)+E\widetilde{g}(\nabla_{E_2}E_2,N)}=0
.
\end{equation}
Let us consider a translation surface $M$ of type I in $\mathbb{R}^3$ parametrized  by $\mathcal{F}(u,v)=(u,v,f(u)+g(v))$. The tangent plane of $M$ is spanned by
\begin{equation}
\mathcal{F}_u=X_1+f'(u)X_3, ~~{\rm and}~~\mathcal{F}_v=X_2+g'(v)X_3,
\end{equation}
while the unit normal $N$ (up to orientation) is given by
\begin{equation}
N=\frac{1}{\alpha}[-f'(u)X_1-g'(v)X_2+X_3],
\end{equation}
where $\alpha^2=f'(u)^2+g'(v)^2+1$.\\
\indent We obtain the coefficients of first fundamental form of $\mathcal{F}$ as
\begin{equation}
E=1+f'(u)^2, ~~~F=f'(u)g'(v),~~~G=1+g'(v)^2.
\end{equation}
Then, the semi-symmetric metric connection (2.1) on the surface is given by
\begin{align}
\left\{\begin{array}{l}
\nabla_{\mathcal{F}_u}\mathcal{F}_u=f'(u)X_1+[f''(u)-1]X_3,\\
\nabla_{\mathcal{F}_u}\mathcal{F}_v=g'(v)X_1,\\
\nabla_{\mathcal{F}_v}\mathcal{F}_u=f'(u)X_2,\\
\nabla_{\mathcal{F}_v}\mathcal{F}_v=g'(v)X_2+[g''(v)-1]X_3.\\
\end{array}\right.
\end{align}
Consequently, the minimality condition (2.9) may be expressed as follows:
\begin{equation}
f''(u)g'(v)^2-2f'(u)^2-2g'(v)^2+f'(u)^2g''(v)+f''(u)+g''(v)-2=0.
\end{equation}
We will solve (2.14). Let us assume first that $f',g',f''$ and $g''$ are different from zero at every point. Taking successive derivatives with respect to $u$ and $v$, we obtain
\begin{equation}
\frac{f'''}{f'f''}=-\frac{g'''}{g'g''}.
\end{equation}
In the following, $c_j$ denotes a constant where $j$ is a positive integer. Remark that the left-hand side of equation (2.15) is a function of $u$, while the right-hand side is a function of $v$. Therefore, there exist three constants
$c_0,c_1,c_2$ such that
\begin{equation}
f''=\frac{c_0}{2}f'^2+c_1,~~~g''=-\frac{c_0}{2}g'^2+c_2.
\end{equation}
Now, plugging (2.16) into (2.14) we obtain
\begin{equation}
(c_2+\frac{c_0}{2}-2)f'^2+c_1+c_2-2=(-c_1+2+\frac{c_0}{2})g'^2.
\end{equation}
Remark that the left-hand side of equation (2.17) is a function of $u$, while the right-hand side is a function of $v$. Therefore,
$(-c_1+2+\frac{c_0}{2})g'^2=\overline{c}_0$. Then $g'^2=\widetilde{c}_0$ or $-c_1+2+\frac{c_0}{2}=0$. But we assume that $g'$ and $g''$ are different from zero at every point
, so we get  \begin{equation}
-c_1+2+\frac{c_0}{2}=0.
\end{equation}
By (2.17), we get
\begin{equation}
(c_2+\frac{c_0}{2}-2)f'^2+c_1+c_2-2=0.
\end{equation}
Since we assume that $f'$ and $f''$ are different from zero at every point, so we get
\begin{equation}
c_2+\frac{c_0}{2}-2=0,~~~c_1+c_2-2=0.
\end{equation}
By (2.18) and (2.20), we have a contradiction. So in this case, we have no solutions.\\
\indent Case 1) There exists $u_0$ such that $f''(u_0)\neq 0$ and there exists $v_0$ such that $g''(v_0)\neq 0$. Then there is an open interval $U$ of $u_0$ such that
$ f''|_U\neq 0$. Then there exists $u_1\in U$ such that $f'(u_1)\neq 0$ and there is an open interval $U_1\subset U$ of $u_1$ such that
$ f'|_{U_1}\neq 0$ and $ f''|_{U_1}\neq 0$. Similarly, there is an open interval $V_1$  such that
$ g'|_{V_1}\neq 0$ and $ g''|_{V_1}\neq 0$. By the above discussions, we know that we have no solutions in this case.\\
\indent Case 2) There exists $u_0$ such that $f''(u_0)\neq 0$ and $g''(v)=0$. So $g'(v)=c_3$ and $g(v)=c_3v+c_4$. By (2.14), we have
\begin{equation}
f''(u)-\frac{2}{c_3^2+1}f'(u)^2-2=0.
\end{equation}
The general solution of this ODE is found as
\begin{equation}
f(u)=-\frac{c_3^2+1}{2}{\rm ln}\mid{\rm cos}(\frac{2}{\sqrt{c_3^2+1}}u-a)\mid+b
\end{equation}
where $a,b$ are constant. So
\begin{equation}
\mathcal{F}(u,v)=(u,v,-\frac{c_3^2+1}{2}{\rm ln}\mid{\rm cos}(\frac{2}{\sqrt{c_3^2+1}}u-a)\mid+c_3v+c_5).
\end{equation}
\indent Case 3) $f''(u)=0$ and there exists $v_0$ such that $g''(v_0)\neq 0$. So $f'(u)=\overline{c}_3$ and $f(u)=\overline{c}_3u+\overline{c}_4$.
Similar to case 2), we have
\begin{equation}
\mathcal{F}(u,v)=(u,v,-\frac{\overline{c_3}^2+1}{2}{\rm ln}\mid{\rm cos}(\frac{2}{\sqrt{\overline{c_3}^2+1}}v-a_1)\mid+\overline{c_3}u+c_6).
\end{equation}
where $\overline{c_3},a_1$ are constant. \\
\indent Case 4) $f''(u)=0$ and  $g''(v)=0$. By (2.14), we have $-2f'(u)^2-2g'(v)^2-2=0$. This is a contradiction, so we have no solutions in this case.\\
 \indent So we have the following theorem£º
\begin{thm}
Type I minimal translation surfaces with respect to $\nabla$ in $\mathbb{R}^3$ are of the forms (2.23) and (2.24).
\end{thm}
\vskip 0.5 true cm
In the following, we obtain all Type II minimal translation surfaces with respect to $\nabla$ in $\mathbb{R}^3$. Let $M$ be a translation surface of type II
parametrized by $\mathcal{F}(u,v)=(u,f(u)+g(v),v)$.
The tangent plane of $M$ is spanned by
\begin{equation}
\mathcal{F}_u=X_1+f'(u)X_2, ~~{\rm and}~~\mathcal{F}_v=g'(v)X_2+X_3,
\end{equation}
while the unit normal $N$ (up to orientation) is given by
\begin{equation}
N=\frac{1}{\beta}[f'(u)X_1-X_2+g'(v)X_3],
\end{equation}
where $\beta^2=f'(u)^2+g'(v)^2+1$.\\
\indent We obtain the coefficients of first fundamental form of $\mathcal{F}$ as
\begin{equation}
E=1+f'(u)^2, ~~~F=f'(u)g'(v),~~~G=1+g'(v)^2.
\end{equation}
Then, the semi-symmetric metric connection (2.1) on the surface is given by
\begin{align}
\left\{\begin{array}{l}
\nabla_{\mathcal{F}_u}\mathcal{F}_u=f''(u)X_2-[f'(u)^2+1]X_3,\\
\nabla_{\mathcal{F}_u}\mathcal{F}_v=X_1+f'(u)X_2-f'(u)g'(v)X_3,\\
\nabla_{\mathcal{F}_v}\mathcal{F}_u=-f'(u)g'(v)X_3,\\
\nabla_{\mathcal{F}_v}\mathcal{F}_v=[g'(v)+g''(v)]X_2-g'(v)^2X_3.\\
\end{array}\right.
\end{align}
Consequently, the minimality condition (2.9) may be expressed as follows:
\begin{equation}
2g'^3+2f'^2g'+g'^2f''+f'^2g''+f''+g''+2g'=0.
\end{equation}
We will solve (2.29). Let us assume first that $f',g',f''$ and $g''$ are different from zero at every point. Taking successive derivatives with respect to $u$ and $v$, we obtain
\begin{equation}
\frac{f'''}{f'f''}=-\frac{2g''+g'''}{g'g''}.
\end{equation}
Remark that the left-hand side of equation (2.30) is a function of $u$, while the right-hand side is a function of $v$. Therefore, there exist two constants
$c_0,c_1$ such that
\begin{equation}
f''=\frac{c_0}{2}f'^2+c_1,~~~\frac{g'''}{g''}=-{c_0}g'-2.
\end{equation}
For (2.29), taking derivative with respect to $v$, then we have
\begin{equation}
6g'^2g''+2f'^2g''+2g'g''f''+f'^2g'''+g'''+2g''=0.
\end{equation}
Plugging (2.31) into (2.32) and $g',g''$ are different from zero , we obtain $g'=-\frac{2c_1-c_0}{6}$. Then $g''=0$. This is a contradiction, so we have no solutions in this case.\\
\indent Case 1) There exists $u_0$ such that $f''(u_0)\neq 0$ and there exists $v_0$ such that $g''(v_0)\neq 0$. By the above discussions, we know that we have no solutions in this case.\\
\indent Case 2) There exists $u_0$ such that $f''(u_0)\neq 0$ and $g''(v)=0$. So $g'(v)=\widetilde{c_0}$ and $g(v)=\widetilde{c_0}v+\widetilde{c_1}$.
 By (2.29), we have
\begin{equation}
f''(u)+\frac{2\widetilde{c_0}}{\widetilde{c_0}^2+1}f'(u)^2+2\widetilde{c_0}=0.
\end{equation}
Since $f''\neq 0$, so $\widetilde{c_0}\neq 0$.
The general solution of this ODE is found as
\begin{equation}
f(u)=\frac{\widetilde{c_0}^2+1}{2\widetilde{c_0}}{\rm ln}\mid{\rm cos}(\frac{2\widetilde{c_0}}{\sqrt{\widetilde{c_0}^2+1}}u-\widetilde{a})\mid+\widetilde{b}
\end{equation}
where $\widetilde{a},\widetilde{b}$ are constant. So
\begin{equation}
\mathcal{F}(u,v)=(u,\frac{\widetilde{c_0}^2+1}{2\widetilde{c_0}}{\rm ln}\mid{\rm cos}(\frac{2\widetilde{c_0}}{\sqrt{\widetilde{c_0}^2+1}}u-\widetilde{a})\mid+\widetilde{c_0}v+\overline{b},v).
\end{equation}
\indent Case 3) $f''(u)=0$ and there exists $v_0$ such that $g''(v_0)\neq 0$. So $f'(u)=\widehat{c}_0$ and $f(u)=\widehat{c}_0u+\widehat{c}_1$.
By (2.29), we have
\begin{equation}
g''+\frac{2}{\widehat{c_0}+1}g'^3+2g'=0.
\end{equation}
Let $g'=h$, then we get $h'=-\frac{2}{\widehat{c_0}+1}h^3-2h$. By $g''(v_0)\neq 0$, we may set $h\neq 0$ on an open interval. Let $W=h^{-2}$, then
 \begin{equation}
W'=\frac{4}{\widehat{c_0}^2+1}+4W.
\end{equation}
Then $W=\widehat{a}e^{4v}-\frac{1}{\widehat{c_0}^2+1}$ where $\widehat{a}>0$. so we get
 \begin{equation}
g(v)=\int_0^v\frac{\pm 1}{\sqrt{\widehat{a}e^{4x}-\frac{1}{\widehat{c_0}^2+1}}}dx+\widehat{b}.
\end{equation}
where $\widehat{a},\widehat{b}$ are constant.
So
\begin{equation}
\mathcal{F}(u,v)=(u,\int_0^v\frac{\pm 1}{\sqrt{\widehat{a}e^{4x}-\frac{1}{\widehat{c_0}^2+1}}}dx+\widehat{c_0}u+b_0,v).
\end{equation}

\indent Case 4) $f''(u)=0$ and  $g''(v)=0$. Then $f'(u)=c_0'$ and  $g'(v)=c_1'$. By (2.29), we get $c_1'=0$. Then
\begin{equation}
\mathcal{F}(u,v)=(u,c_0'u+b',v).
\end{equation}
where $c_0',b'$ are constant.\\

 \indent So we have the following theorem£º
\begin{thm}
Type II minimal translation surfaces with respect to $\nabla$ in $\mathbb{R}^3$ are of the forms (2.35),(2.39) and (2.40).
\end{thm}
\vskip 0.5 true cm
In the following, we consider Type III minimal translation surfaces with respect to $\nabla$ in $\mathbb{R}^3$. Let $M$ be a translation surface of type III
parametrized by $\mathcal{F}(u,v)=(f(u)+g(v),u,v)$.
The tangent plane of $M$ is spanned by
\begin{equation}
\mathcal{F}_u=f'(u)X_1+X_2, ~~{\rm and}~~\mathcal{F}_v=g'(v)X_1+X_3,
\end{equation}
while the unit normal $N$ (up to orientation) is given by
\begin{equation}
N=\frac{1}{\gamma}[X_1-f'(u)X_2-g'(v)X_3],
\end{equation}
where $\gamma^2=f'(u)^2+g'(v)^2+1$.\\
\indent We obtain the coefficients of first fundamental form of $\mathcal{F}$ as
\begin{equation}
E=1+f'(u)^2, ~~~F=f'(u)g'(v),~~~G=1+g'(v)^2.
\end{equation}
Then, the semi-symmetric metric connection (2.1) on the surface is given by
\begin{align}
\left\{\begin{array}{l}
\nabla_{\mathcal{F}_u}\mathcal{F}_u=f''(u)X_1-[f'(u)^2+1]X_3,\\
\nabla_{\mathcal{F}_u}\mathcal{F}_v=f'(u)X_1+X_2-f'(u)g'(v)X_3,\\
\nabla_{\mathcal{F}_v}\mathcal{F}_u=-f'(u)g'(v)X_3,\\
\nabla_{\mathcal{F}_v}\mathcal{F}_v=[g'(v)+g''(v)]X_1-g'(v)^2X_3.\\
\end{array}\right.
\end{align}
Consequently, the minimality condition (2.9) may be expressed as follows:
\begin{equation}
2g'^3+2f'^2g'+g'^2f''+f'^2g''+f''+g''+2g'=0.
\end{equation}
(2.45) is the same as (2.29), so similar to Theorem 2.3, we can get Type III minimal translation surfaces with respect to $\nabla$ in $\mathbb{R}^3$.\\
\indent We define a special semi-symmetric non-metric connection by
\begin{equation}
\overline{\nabla}_{X}Y=\nabla^L_{X}Y+\widetilde{g}(Y,X_3)X.
\end{equation}
The connection $\overline{\nabla}$ of $\mathbb{R}^3$ is given by
\begin{align}
&\overline{\nabla}_{X_1}X_1=0,~~\overline{\nabla}_{X_1}X_2=0,~~\overline{\nabla}_{X_1}X_3=X_1,\\\notag
&\overline{\nabla}_{X_2}X_1=0,~~\overline{\nabla}_{X_2}X_2=0,~~\overline{\nabla}_{X_2}X_3=X_2,\\\notag
&\overline{\nabla}_{X_3}X_1=0,~~\overline{\nabla}_{X_3}X_2=0,~~\overline{\nabla}_{X_3}X_3=X_3.
\notag
\end{align}

In the following, we consider type I minimal translation surfaces with respect to $\overline{\nabla}$ in $\mathbb{R}^3$.
For $\mathcal{F}(u,v)=(u,v,f(u)+g(v))$. then $\mathcal{F}_u, \mathcal{F}_v,N,E,F,G$ is computed by (2.10)-(2.12).

Then, the semi-symmetric non-metric connection (2.46) on the surface is given by
\begin{align}
\left\{\begin{array}{l}
\overline{\nabla}_{\mathcal{F}_u}\mathcal{F}_u=f'(u)X_1+[f''(u)+f'(u)^2]X_3,\\
\overline{\nabla}_{\mathcal{F}_u}\mathcal{F}_v=g'(v)X_1+f'(u)g'(v)X_3,\\
\overline{\nabla}_{\mathcal{F}_v}\mathcal{F}_u=f'(u)X_2+f'(u)g'(v)X_3,\\
\overline{\nabla}_{\mathcal{F}_v}\mathcal{F}_v=g'(v)X_2+[g'(v)^2+g''(v)]X_3.\\
\end{array}\right.
\end{align}
Similar to (2.9), we have the minimality condition with respect to $\overline{\nabla}$.
 Consequently, the minimality condition may be expressed as follows:
\begin{equation}
\frac{f''(u)}{1+f'(u)^2}=-\frac{g''(v)}{1+g'(v)^2}.
\end{equation}
Solving (2.49), we obtain

\begin{thm}
Type I minimal translation surfaces with respect to $\overline{\nabla}$ in $\mathbb{R}^3$ are of the following forms 
\begin{equation}
\mathcal{F}(u,v)=(u,v,c_0u+c_1v+c_2).
\end{equation}
\begin{equation}
\mathcal{F}(u,v)=(u,v,\frac{1}{c}{\rm ln}\frac{\mid{\rm cos}(cu-c_3)\mid}{\mid{\rm cos}(cv-c_4)\mid}+c_5),
\end{equation}
where $c\neq 0$.
\end{thm}
\vskip 0.5 true cm
For Type II and III minimal translation surfaces with respect to $\overline{\nabla}$ in $\mathbb{R}^3$, we also get the equation (2.49) and we have Theorems
similar to Theorem 2.4.

\section{Minimal translation surfaces with respect to semi-symmetric connections in $\mathbb{R}_1^3$ }
\indent Let $\mathbb{R}_1^3$ be the $3$-dimensional Minkowski space with the canonical Minkowski metric $\widetilde{g_1}=dx^2+dy^2-dz^2$.  Let $\nabla^L_{X_i}X_j=0$ for $1\leq i,j\leq 3$ be the Levi-Civita connection on
$\mathbb{R}_1^3$. We define a special semi-symmetric metric connection by
\begin{equation}
\nabla_{X}Y=\nabla^L_{X}Y+\widetilde{g_1}(Y,X_3)X-\widetilde{g_1}(X,Y)X_3.
\end{equation}
The connection $\nabla$ of $\mathbb{R}_1^3$ is given by
\begin{align}
&\nabla_{X_1}X_1=-X_3,~~\nabla_{X_1}X_2=0,~~\nabla_{X_1}X_3=-X_1,\\\notag
&\nabla_{X_2}X_1=0,~~\nabla_{X_2}X_2=-X_3,~~\nabla_{X_2}X_3=-X_2,\\\notag
&\nabla_{X_3}X_1=0,~~\nabla_{X_3}X_2=0,~~\nabla_{X_3}X_3=0.
\notag
\end{align}
In the following, we assume that $M$ is a spacelike surface of $\mathbb{R}_1^3$, that is the induced metric on $M$ is Riemannian metric. When 
 $M$ is a timelike surface of $\mathbb{R}_1^3$, we have similar discussions.
Let us consider a translation surface $M$ of type I in $\mathbb{R}_1^3$ parametrized  by $\mathcal{F}(u,v)=(u,v,f(u)+g(v))$. 
$\mathcal{F}_u$ and $\mathcal{F}_v$ are given by (2.10).
While the unit normal $N$ (up to orientation) is given by
\begin{equation}
N=\frac{1}{\alpha_1}[-f'(u)X_1-g'(v)X_2-X_3],
\end{equation}
where $\alpha_1^2=1-f'(u)^2-g'(v)^2$.\\
\indent We obtain the coefficients of first fundamental form of $\mathcal{F}$ as
\begin{equation}
E=1-f'(u)^2, ~~~F=-f'(u)g'(v),~~~G=1-g'(v)^2.
\end{equation}
Then, the semi-symmetric metric connection (3.1) on the surface is given by
\begin{align}
\left\{\begin{array}{l}
\nabla_{\mathcal{F}_u}\mathcal{F}_u=-f'(u)X_1+[f''(u)-1]X_3,\\
\nabla_{\mathcal{F}_u}\mathcal{F}_v=-g'(v)X_1,\\
\nabla_{\mathcal{F}_v}\mathcal{F}_u=-f'(u)X_2,\\
\nabla_{\mathcal{F}_v}\mathcal{F}_v=-g'(v)X_2+[g''(v)-1]X_3.\\
\end{array}\right.
\end{align}
In this case, we have the same minimality condition with (2.9). So the minimality condition may be expressed as follows:
\begin{equation}
f''(u)g'(v)^2-2f'(u)^2-2g'(v)^2+f'(u)^2g''(v)-f''(u)-g''(v)+2=0.
\end{equation}
Similar to the discussions of the case 1) in the page 4, we get\\
 \indent Case 1) There exists $u_0$ such that $f''(u_0)\neq 0$ and there exists $v_0$ such that $g''(v_0)\neq 0$. We know that we have no solutions in this case.\\
\indent Case 2) There exists $u_0$ such that $f''(u_0)\neq 0$ and $g''(v)=0$. So $g'(v)=c$ and $g(v)=cv+\overline{c}$.
By (3.6), we have
\begin{equation}
(c^2-1)f''(u)-2f'(u)^2-2(c^2-1)=0.
\end{equation}
If $c^2=1$, then $f'=0$ and $f''=0$. This is a contradiction. So $c^2\neq 1$ and we get
\begin{equation}
f''(u)-\frac{2}{c^2-1}f'(u)^2-2=0.
\end{equation}
If $c^2>1$, then
the general solution of this ODE (3.8) is found as
\begin{equation}
f(u)=-\frac{c^2-1}{2}{\rm ln}\mid{\rm cos}(\frac{2}{\sqrt{c^2-1}}u-a)\mid+b
\end{equation}
where $a,b$ are constant. So
\begin{equation}
\mathcal{F}(u,v)=(u,v,-\frac{c^2-1}{2}{\rm ln}\mid{\rm cos}(\frac{2}{\sqrt{c^2-1}}u-a)\mid+cv+\overline{b}).
\end{equation}
If $c^2<1$, then
the general solution of this ODE (3.8) is found as
\begin{equation}
f(u)=\int_0^u\frac{\sqrt{1-c^2}+\widetilde{c}e^{-\frac{4}{\sqrt{1-c^2}}x}}{1-\widetilde{c}e^{-\frac{4}{\sqrt{1-c^2}}x}}dx
\end{equation}
where $\widetilde{c}$ is a nonzero constant. 
So
\begin{equation}
\mathcal{F}(u,v)=(u,v,\int_0^u\frac{\sqrt{1-c^2}+\widetilde{c}e^{-\frac{4}{\sqrt{1-c^2}}x}}{1-\widetilde{c}e^{-\frac{4}{\sqrt{1-c^2}}x}}dx+cv+\widetilde{b}).
\end{equation}
\indent Case 3) $f''(u)=0$ and there exists $v_0$ such that $g''(v_0)\neq 0$. So $f'(u)=\widehat{c}$.
Similar to case 2), we have if $\widehat{c}^2>1$, then
\begin{equation}
\mathcal{F}(u,v)=(u,v,-\frac{\widehat{c}^2-1}{2}{\rm ln}\mid{\rm cos}(\frac{2}{\sqrt{\widehat{c}^2-1}}v-a_1)\mid+\widehat{c}u+\overline{b}_1).
\end{equation}
If $\widehat{c}^2<1$, then
\begin{equation}
\mathcal{F}(u,v)=(u,v,\int_0^v\frac{\sqrt{1-\widehat{c}^2}+\widetilde{c}_1e^{-\frac{4}{\sqrt{1-\widehat{c}^2}}x}}{1-\widetilde{c}_1e^{-\frac{4}
{\sqrt{1-\widehat{c}^2}}x}}dx+\widehat{c}u+\widetilde{b}).
\end{equation}
\indent Case 4) $f''(u)=0$ and  $g''(v)=0$. By (3.6), we have $1-f'(u)^2-g'(v)^2=0$. By $N$ is timelike, then $1-f'(u)^2-g'(v)^2>0$. This is a contradiction, so we have no solutions in this case.\\
 \indent So we have the following theorem£º
\begin{thm}
Type I minimal translation surfaces with respect to $\nabla$ in $\mathbb{R}_1^3$ are of the forms (3.10),(3.12),(3.13) and (3.14).
\end{thm}
\vskip 0.5 true cm
 \indent In the following, we obtain all Type II minimal translation surfaces with respect to $\nabla$ in $\mathbb{R}_1^3$. Let $M$ be a translation surface of type II
parametrized by $\mathcal{F}(u,v)=(u,f(u)+g(v),v)$. $\mathcal{F}_u$ and $\mathcal{F}_v$ are given by (2.25).
While the unit normal $N$ (up to orientation) is given by
\begin{equation}
N=\frac{1}{\beta_1}[f'(u)X_1-X_2-g'(v)X_3],
\end{equation}
where $\beta_1^2=-1-f'(u)^2+g'(v)^2$.\\
\indent We obtain the coefficients of first fundamental form of $\mathcal{F}$ as
\begin{equation}
E=1+f'(u)^2, ~~~F=f'(u)g'(v),~~~G=g'(v)^2-1.
\end{equation}
Then, the semi-symmetric metric connection (3.1) on the surface is given by
\begin{align}
\left\{\begin{array}{l}
\nabla_{\mathcal{F}_u}\mathcal{F}_u=f''(u)X_2-[f'(u)^2+1]X_3,\\
\nabla_{\mathcal{F}_u}\mathcal{F}_v=-X_1-f'(u)X_2-f'(u)g'(v)X_3,\\
\nabla_{\mathcal{F}_v}\mathcal{F}_u=-f'(u)g'(v)X_3,\\
\nabla_{\mathcal{F}_v}\mathcal{F}_v=[-g'(v)+g''(v)]X_2-g'(v)^2X_3.\\
\end{array}\right.
\end{align}
Consequently, the minimality condition may be expressed as follows:
\begin{equation}
2g'^3-2f'^2g'+g'^2f''+f'^2g''-f''+g''-2g'=0.
\end{equation}
We will solve (3.18). Let us assume first that $f',g',f''$ and $g''$ are different from zero at every point. Taking successive derivatives with respect to $u$ and $v$, we obtain
\begin{equation}
\frac{f'''}{f'f''}=\frac{-2g''+g'''}{g'g''}.
\end{equation}
 Therefore, there exist two constants
$c_0,c_1$ such that
\begin{equation}
f''=\frac{c_0}{2}f'^2+c_1,~~~\frac{g'''}{g''}={c_0}g'+2.
\end{equation}
For (3.18), taking derivative with respect to $v$, then we have
\begin{equation}
6g'^2g''-2f'^2g''+2g'g''f''+f'^2g'''+g'''-2g''=0.
\end{equation}
Plugging (3.20) into (3.21) and $g',g''$ are different from zero , we obtain
$6g'=-2c_0f'^2-2c_1-c_0$. Then $g''=0$. This is a contradiction, so we have no solutions in this case.\\
\indent Case 1) There exists $u_0$ such that $f''(u_0)\neq 0$ and there exists $v_0$ such that $g''(v_0)\neq 0$. By the above discussions, we know that we have no solutions in this case.\\
\indent Case 2) There exists $u_0$ such that $f''(u_0)\neq 0$ and $g''(v)=0$. So $g'(v)=\widetilde{c_0}$ and $g(v)=\widetilde{c_0}v+\widetilde{c_1}$.
  By (3.18), we have
\begin{equation}
(\widetilde{c_0}^2-1)f''(u)-2\widetilde{c_0}f'(u)^2+2\widetilde{c_0}(\widetilde{c_0}^2-1)=0.
\end{equation}
Since $f''\neq 0$, so $\widetilde{c_0}\neq 0$ and $\widetilde{c_0}^2\neq 1$. So
\begin{equation}
f''(u)-\frac{2\widetilde{c_0}}{\widetilde{c_0}^2-1}f'(u)^2+2\widetilde{c_0}=0.
\end{equation}
The general solution of this ODE (3.23) is found as if $\widetilde{c_0}^2< 1$
\begin{equation}
f(u)=\frac{1-\widetilde{c_0}^2}{2\widetilde{c_0}}{\rm ln}\mid{\rm cos}(\frac{2\widetilde{c_0}}{\sqrt{1-\widetilde{c_0}^2}}u-\widetilde{a})\mid+\widetilde{b}
\end{equation}
where $\widetilde{a},\widetilde{b}$ are constant. So
\begin{equation}
\mathcal{F}(u,v)=(u,\frac{1-\widetilde{c_0}^2}{2\widetilde{c_0}}{\rm ln}\mid{\rm cos}(\frac{2\widetilde{c_0}}{\sqrt{1-\widetilde{c_0}^2}}u-\widetilde{a})\mid+\widetilde{c_0}v+\overline{b},v).
\end{equation}
If $\widetilde{c_0}^2>1$
then
the general solution of this ODE (3.23) is found as
\begin{equation}
f(u)=\int_0^u\frac{\sqrt{\widetilde{c_0}^2-1}+{c_1}e^{\frac{4\widetilde{c_0}}{\sqrt{\widetilde{c_0}^2-1}}x}}{1-{c_1}e^{\frac{4\widetilde{c_0}}{\sqrt{
\widetilde{c_0}^2-1}}x}}dx,
\end{equation}
where ${c_1}$ is a nonzero constant. 
So
\begin{equation}
\mathcal{F}(u,v)=(u,\int_0^u\frac{\sqrt{\widetilde{c_0}^2-1}+{c_1}e^{\frac{4\widetilde{c_0}}{\sqrt{\widetilde{c_0}^2-1}}x}}{1-{c_1}e^{\frac{4\widetilde{c_0}}{\sqrt{
\widetilde{c_0}^2-1}}x}}dx+\widetilde{c_0}v+\overline{b_1},v).
\end{equation}
\indent Case 3) $f''(u)=0$ and there exists $v_0$ such that $g''(v_0)\neq 0$. So $f'(u)=\widehat{c}_0$ and $f(u)=\widehat{c}_0u+\widehat{c}_1$.
By (3.18), we have
\begin{equation}
g''+\frac{2}{\widehat{c_0}^2+1}g'^3-2g'=0.
\end{equation}
We get
 \begin{equation}
g(v)=\int_0^v\frac{\pm 1}{\sqrt{\widehat{a}e^{-4x}+\frac{1}{\widehat{c_0}^2+1}}}dx+\widehat{b}.
\end{equation}
where $\widehat{a}\neq 0,\widehat{b}$ are constant.
So
\begin{equation}
\mathcal{F}(u,v)=(u,\int_0^v\frac{\pm 1}{\sqrt{\widehat{a}e^{-4x}+\frac{1}{\widehat{c_0}^2+1}}}dx+\widehat{c_0}u+b_0,v).
\end{equation}
\indent Case 4) $f''(u)=0$ and  $g''(v)=0$. Then $f'(u)=c_0'$ and  $g''(v)=c_1'$. By (3.18), we get $c_1'(c_1'^2-c_0'^2-2)=0$. Then
\begin{equation}
\mathcal{F}(u,v)=(u,c_0'u+c_1'v+b',v).
\end{equation}
where $c_1'=0$ or $(c_1'^2-c_0'^2-2)=0.$\\
\indent So we have the following theorem£º
\begin{thm}
Type II minimal translation surfaces with respect to $\nabla$ in $\mathbb{R}_1^3$ are of the forms (3.25),(3.27).(3.30) and (3.31).
\end{thm}
\vskip 0.5 true cm
For Type III minimal translation surfaces with respect to $\nabla$ in $\mathbb{R}_1^3$, we also get (3.18) and we have Theorem similar to Theorem 3.2.\\
\indent We define a special semi-symmetric non-metric connection in $\mathbb{R}_1^3$ by
\begin{equation}
\overline{\nabla}_{X}Y=\nabla^L_{X}Y+\widetilde{g_1}(Y,X_3)X.
\end{equation}
The connection $\overline{\nabla}$ of $\mathbb{R}_1^3$ is given by
\begin{align}
&\overline{\nabla}_{X_1}X_1=0,~~\overline{\nabla}_{X_1}X_2=0,~~\overline{\nabla}_{X_1}X_3=-X_1,\\\notag
&\overline{\nabla}_{X_2}X_1=0,~~\overline{\nabla}_{X_2}X_2=0,~~\overline{\nabla}_{X_2}X_3=-X_2,\\\notag
&\overline{\nabla}_{X_3}X_1=0,~~\overline{\nabla}_{X_3}X_2=0,~~\overline{\nabla}_{X_3}X_3=-X_3.
\notag
\end{align}
In the following, we consider type I minimal translation surfaces with respect to $\overline{\nabla}$ in $\mathbb{R}_1^3$.
For $\mathcal{F}(u,v)=(u,v,f(u)+g(v))$, then $\mathcal{F}_u, \mathcal{F}_v,N,E,F,G$ is computed by (3.3) and (3.4).
Then, the semi-symmetric non-metric connection (3.32) on the surface is given by
\begin{align}
\left\{\begin{array}{l}
\overline{\nabla}_{\mathcal{F}_u}\mathcal{F}_u=-f'(u)X_1+[f''(u)-f'(u)^2]X_3,\\
\overline{\nabla}_{\mathcal{F}_u}\mathcal{F}_v=-g'(v)X_1-f'(u)g'(v)X_3,\\
\overline{\nabla}_{\mathcal{F}_v}\mathcal{F}_u=-f'(u)X_2-f'(u)g'(v)X_3,\\
\overline{\nabla}_{\mathcal{F}_v}\mathcal{F}_v=-g'(v)X_2+[-g'(v)^2+g''(v)]X_3.\\
\end{array}\right.
\end{align}
Similar to (2.9), we have the minimality condition with respect to $\overline{\nabla}$ in $\mathbb{R}_1^3$.
$N$ is timelike, so $1-f'(u)^2-g'(v)^2>0$ and $1-f'(u)^2>0$ and $1-g'(v)^2>0$. Consequently, the minimality condition may be expressed as follows:
\begin{equation}
\frac{f''(u)}{1-f'(u)^2}=-\frac{g''(v)}{1-g'(v)^2}=c_0.
\end{equation}
When $c_0=0$, we get
\begin{equation}
\mathcal{F}(u,v)=(u,v,c_1u+c_2v+c_3).
\end{equation}
When $c_0\neq 0$, we get
\begin{equation}
{f''(u)}+c_0f'(u)^2-c_0=0,~~~{g''(u)}-c_0g'(u)^2+c_0=0.
\end{equation}
Solving (3.37), we obtain
\begin{thm}
Type I minimal translation surfaces with respect to $\overline{\nabla}$ in $\mathbb{R}_1^3$ are (3.36) and the following form
\begin{equation}
\mathcal{F}(u,v)=(u,v,\frac{1}{c_0}{\rm ln}\frac{\mid e^{c_0u}-\widehat{c}e^{-c_0u}\mid}{\mid e^{-c_0v}-\widehat{c_1}e^{c_0v} \mid}+a),
\end{equation}
where $c_0,\widehat{c},\widehat{c_1}\neq 0$.
\end{thm}
\vskip 0.5 true cm
\indent For Type II minimal translation surfaces with respect to $\overline{\nabla}$ in $\mathbb{R}_1^3$,
and $\mathcal{F}(u,v)=(u,f(u)+g(v),v)$, then $\mathcal{F}_u, \mathcal{F}_v,N,E,F,G$ is computed by (3.15) and (3.16).
Then, the semi-symmetric non-metric connection (3.32) on the surface is given by
\begin{align}
\left\{\begin{array}{l}
\overline{\nabla}_{\mathcal{F}_u}\mathcal{F}_u=f''(u)X_2,\\
\overline{\nabla}_{\mathcal{F}_u}\mathcal{F}_v=-X_1-f'(u)X_2,\\
\overline{\nabla}_{\mathcal{F}_v}\mathcal{F}_u=0,\\
\overline{\nabla}_{\mathcal{F}_v}\mathcal{F}_v=[-g'(v)+g''(v)]X_2-X_3.\\
\end{array}\right.
\end{align}
 Consequently, the minimality condition may be expressed as follows:
\begin{equation}
\frac{f''(u)}{1+f'(u)^2}=\frac{g''(v)}{1-g'(v)^2}=\overline{c_0}.
\end{equation}
When $\overline{c_0}=0$, we get
\begin{equation}
\mathcal{F}(u,v)=(u,c_1u+c_2v+c_3,v).
\end{equation}
When $\overline{c_0}\neq 0$, we get
\begin{equation}
{f''(u)}-\overline{c_0}f'(u)^2-\overline{c_0}=0,~~~{g''(u)}+\overline{c_0}g'(u)^2-\overline{c_0}=0.
\end{equation}
Solving (3.42), we obtain
\begin{thm}
Type II minimal translation surfaces with respect to $\overline{\nabla}$ in $\mathbb{R}_1^3$ are (3.41) and the following form
\begin{equation}
\mathcal{F}(u,v)=(u,\frac{1}{\overline{c_0}}{\rm ln}\frac{\mid e^{\overline{c_0}v}-c_3e^{-\overline{c_0}v}\mid}{\mid
{\rm cos}(\overline{c_0}u+c_4)
\mid}+b,v),
\end{equation}
where $\overline{c_0},c_3\neq 0$.
\end{thm}
\vskip 0.5 true cm
For Type III minimal translation surfaces with respect to $\overline{\nabla}$ in $\mathbb{R}_1^3$, we also get (3.40) and we have Theorem similar to Theorem 3.4.\\
\vskip 1 true cm

\section{Acknowledgements}

The author was supported in part by  NSFC No.11771070. The author thanks Dr. Sining Wei for her helpful discussions.

\vskip 1 true cm


\bigskip
\bigskip

\noindent {\footnotesize {\it Y. Wang} \\
{School of Mathematics and Statistics, Northeast Normal University, Changchun 130024, China}\\
{Email: wangy581@nenu.edu.cn}

\end{document}